\documentclass[10pt,twoside,reqno]{amsart}
\usepackage{amssymb}
\usepackage{multirow}
\calclayout

\numberwithin{equation}{section}
\makeatletter

\renewcommand{\@secnumfont}{\bfseries}

\renewcommand{\section}{\@startsection{section}{1}%
  {0mm}{.7\linespacing\@plus\linespacing}{.5\linespacing}
  {\normalfont\bfseries\centering}}

\newcommand{\bibsection}{\@startsection{section}{1}%
  {0mm}{.7\linespacing\@plus\linespacing}{.5\linespacing}
  {\normalfont\scshape\centering}}

\renewcommand{\@biblabel}[1]{#1.}

\newtheorem{thm}{\bf Theorem}[section]

\begin{document}

\vspace{1.3cm}

\title {$\lambda$-analogues of $r$-stirling numbers of the first kind}

\author{Taekyun Kim}
\address{Department of Mathematics, Kwangwoon University, Seoul 139-701, Republic of Korea}
\email{tkkim@kw.ac.kr}

\author{Dae San Kim}
\address{Department of Mathematics, Sogang University, Seoul 121-742, Republic of Korea}
\email{dskim@sogang.ac.kr}

\subjclass[2010]{11B73; 11B83}
\keywords{$\lambda$-analogues of the $r$-Stirling numbers of the first kind,  higher-order Daehee polynomials}
\begin{abstract} 
In this paper, we study $\lambda$-analogues of the $r$-Stirling numbers of the first kind which have close connections with the $r$-Stirling numbers of the first kind and $\lambda $-Stirling numbers of the first kind. Specifically, we give the recurrence relations for these numbers and show their connections with the $\lambda $-Stirling numbers of the first kind and higher-order Daehee polynomials.
\end{abstract}

\maketitle

\section{Introduction}

It is known that the Stirling numbers of the first kind are defined as
\begin{equation}\begin{split}\label{01}
(x)_n = \sum_{l=0}^n S_1(n,l)x^l,\quad (\text{see}\,\, [1,2,6-9,14]),
\end{split}\end{equation}

\noindent where $(x)_0=1$, $(x)_n=x(x-1)\cdots(x-n+1),\,\,(n \geq 1)$. 

For $\lambda \in \mathbb{R}$, the $\lambda $-analogue of falling factorial sequence is defined by
\begin{equation}\begin{split}\label{02}
&(x)_{0,\lambda }=1, (x)_{n,\lambda }=x(x-\lambda )(x-2\lambda )\cdots(x-(n-1)\lambda ),\,\,(n \geq 1), \\
&(\text{see}\,\, [2,10,14,15,17]).
\end{split}\end{equation}

In view of \eqref{01}, we define $\lambda $-analogues of the Stirling numbers of the first kind as
	
\begin{equation}\begin{split}\label{03}
(x)_{n,\lambda }= \sum_{k=0}^n S_{1,\lambda }(n,k) x^k,\quad (\text{see}\,\, [2,11-13,16,17]).
\end{split}\end{equation}

It is not difficult to show that
\begin{equation}\begin{split}\label{04}
(1+\lambda t)^{\frac{x}{\lambda }} = \sum_{l=0}^\infty {x \choose l}_\lambda  t^l = \sum_{l=0}^\infty \frac{(x)_{l,\lambda }}{l!} t^l,\quad (\text{see}\,\, [4,7-17]),
\end{split}\end{equation}

\noindent where ${x \choose l}_\lambda $ are the $\lambda$-analogues of binomial coefficients ${x \choose n}$ given by ${x \choose l}_\lambda  = \frac{(x)_{l,\lambda }}{l!}$.

The $r$-Stirling numbers of the first kind are defined by the generating function 
\begin{equation}\begin{split}\label{05}
\frac{1}{k!} \big( \log(1+t) \big)^k (1+t)^r = \sum_{n=k}^\infty S_1^{(r)}(n,k) \frac{t^n}{n!},\quad (\text{see}\,\, [3,20-23]).
\end{split}\end{equation}

\noindent where $k \in \mathbb{N}\cup \{0\}$ and $r \in \mathbb{R}$. 

The unsigned $r$-Stirling numbers of the first kind are defined as
\begin{equation}\begin{split}\label{06}
(x+r)(x+r+1)\cdots(x+r+n-1) = \sum_{k=0}^n \genfrac[]{0pt}{1}{n+r}{k+r}_r x^k,\quad (\text{see}\,\, [1,17,22]).
\end{split}\end{equation}

Thus, by \eqref{05}, we get
\begin{equation}\begin{split}\label{07}
(x+r)_n = (x+r)(x+r-1)\cdots(x+r-n+1) = \sum_{k=0}^n S_1^{(r)}(n,k)x^k,\quad (\text{see}\,\, [1]).
\end{split}\end{equation}

From \eqref{05} and \eqref{07}, we note that
\begin{equation}\begin{split}\label{08}
S_1^{(-r)}(n,k) = (-1)^{n-k} \genfrac[]{0pt}{1}{n+r}{k+r}_r.
\end{split}\end{equation}

The higher-order Daehee polynomials are defined by 
\begin{equation}\begin{split}\label{09}
\left( \frac{\log(1+t)}{t} \right)^k (1+t)^x = \sum_{n=0}^\infty   D_n^{(k)}(x) \frac{t^n}{n!},\quad (\text{see}\,\, [5,18,19,24]).
\end{split}\end{equation}

\noindent When $x=0$, $D_n^{(k)} = D_n^{(k)}(0)$ are called the higher-order Daehee numbers. In particular, for $k=1$, $D_n(x) = D_n^{(1)}(x)$, $(n \geq 0)$, are called the ordinary Daehee polynomials.

In this paper, we consider $\lambda $-analogues of $r$-Stirling numbers of the first kind which are derived from the $\lambda $-analogues of the falling factorial sequence and investigate some properties for these numbers. Specifically, we give some identities and recurrence relations for the $\lambda $-analogues of $r$-Stirling numbers of the first kind and show their connections with the $\lambda $-Stirling numbers of the first kind and higher-order Daehee polynomials.

\section{$\lambda $-analogues of $r$-Stirling numbers of the first kind}

From \eqref{03} and \eqref{04}, we have
\begin{equation}\begin{split}\label{10}
(1+\lambda t)^{\frac{x}{\lambda }}&= \sum_{k=0}^\infty (x)_{k,\lambda}\frac{t^k}{k!} 
= \sum_{k=0}^\infty \left( \sum_{n=0}^k S_{1,\lambda }(k,n) x^n \right) \frac{t^k}{k!}\\
& = \sum_{n=0}^\infty \left( n! \sum_{k=n}^\infty S_{1,\lambda }(k,n) \frac{t^k}{k!} \right) \frac{x^n}{n!}.
\end{split}\end{equation}

On the other hand, we also have
\begin{equation}\begin{split}\label{11}
(1+\lambda t)^{\frac{x}{\lambda }} = e^{\frac{x}{\lambda}  \log(1+\lambda t)} = \sum_{n=0}^\infty \left( \frac{\log(1+\lambda t)}{\lambda } \right)^n \frac{x^n}{n!}.
\end{split}\end{equation}

Therefore, by \eqref{10} and \eqref{11}, we get the generating function for $S_{1,\lambda }(n,k)$, $(n,k \geq 0)$, which is given by
\begin{equation}\begin{split}\label{12}
\frac{1}{n!} \left( \frac{\log(1+\lambda t)}{\lambda } \right)^n = \sum_{k=n}^\infty S_{1,\lambda }(k,n) \frac{t^k}{k!}.
\end{split}\end{equation}

Now, we define $\lambda $-analogues of $r$-Stirling numbers of the first kind as
\begin{equation}\begin{split}\label{13}
\frac{1}{k!}  \left( \frac{\log(1+\lambda t)}{\lambda } \right)^k (1+\lambda t)^{\frac{r}{\lambda }} = \sum_{n=k}^\infty S_{1,\lambda }^{(r)}(n,k) \frac{t^n}{n!},
\end{split}\end{equation}

\noindent where $k \in \mathbb{N}\cup \{0\}$, and $r \in \mathbb{R}$.

From \eqref{12} and \eqref{13}, we note that $S_{1,\lambda }^{(0)}(n,k) = S_{1,\lambda }(n,k)$, $(n \geq k \geq 0)$. Also, it  is easy to show that
\begin{equation}\begin{split}\label{14}
(1+\lambda t)^{\frac{x}{\lambda }} (1+\lambda t)^{\frac{r}{\lambda }} = \sum_{n=0}^\infty (x+r)_{n,\lambda } \frac{t^n}{n!}.
\end{split}\end{equation}

By \eqref{14}, we get
\begin{equation}\begin{split}\label{15}
\sum_{n=0}^\infty  (x+r)_{n,\lambda } \frac{t^n}{n!}&= \sum_{n=0}^\infty {x+r \choose n}_\lambda  t^n = (1+\lambda t)^{\frac{r}{\lambda }} e^{\frac{x}{\lambda } \log(1+\lambda t)}\\
&=\sum_{k=0}^{\infty} x^k \frac{1}{k!} \left( \frac{\log(1+\lambda t}{\lambda } \right)^k (1+\lambda t)^{\frac{r}{\lambda }} \\
&= \sum_{k=0}^\infty x^k \sum_{n=k}^\infty S_{1,\lambda }^{(r)}(n,k) \frac{t^n}{n!} = \sum_{n=0}^\infty \left( \sum_{k=0}^n S_{1,\lambda }^{(r)}(n,k)  x^k \right) \frac{t^n}{n!}.
\end{split}\end{equation}

Therefore, by comparing the coefficients on both sides of \eqref{15}, we obtain the following theorem.

\begin{thm}
For $n \geq 0$, we have
\begin{equation*}\begin{split}
(x+r)_{n,\lambda } = \sum_{k=0}^n S_{1,\lambda }^{(r)}(n,k)  x^k.
\end{split}\end{equation*}
\end{thm}

Now, we observe that
\begin{equation}\begin{split}\label{16}
&\sum_{k=0}^\infty x^k \frac{1}{k!} \left( \frac{\log(1+\lambda t)}{\lambda } \right)^k (1+\lambda t)^{\frac{r}{\lambda }}\\
& = \left( \sum_{k=0}^\infty x^k \sum_{m=k}^\infty S_{1,\lambda }(m,k) \frac{t^m}{m!} \right) \left( \sum_{l=0}^\infty (r)_{l,\lambda } \frac{t^l}{l!} \right)\\
&= \left( \sum_{m=0}^\infty \sum_{k=0}^m S_{1,\lambda }(m,k) x^k  \frac{t^m}{m!} \right) \left( \sum_{l=0}^\infty (r)_{l,\lambda } \frac{t^l}{l!} \right) \\&= \sum_{n=0}^\infty \left( \sum_{m=0}^n \sum_{k=0}^m {n \choose m} S_{1,\lambda }(m,k) (r)_{n-m,\lambda } x^k \right) \frac{t^n}{n!}\\
&=\sum_{n=0}^\infty \left( \sum_{k=0}^n \sum_{m=k}^n {n \choose m} S_{1,\lambda }(m,k) (r)_{n-m,\lambda } x^k \right) \frac{t^n }{n!}.
\end{split}\end{equation}
Thus, by \eqref{15} and \eqref{16}, we get
\begin{equation}\begin{split}\label{17}
 \sum_{k=0}^n S_{1,\lambda }^{(r)}(n,k) x^k =  \sum_{k=0}^n \left( \sum_{m=k}^n {n \choose m} S_{1,\lambda }(m,k) (r)_{n-m,\lambda } \right) x^k.
 \end{split}\end{equation} 

Therefore, by comparing the coefficients on both sides of \eqref{17}, we obtain the following theorem.

\begin{thm}
For $n \geq 0$, we have
\begin{equation*}\begin{split}
S_{1,\lambda }^{(r)}(n,k) = \sum_{m=k}^n {n \choose m} S_{1,\lambda }(m,k) (r)_{n-m,\lambda }.
\end{split}\end{equation*}
\end{thm}

Now, we define $\lambda $-analogues of the unsigned $r$-Stirling numbers of the first kind as follows:

\begin{equation}\begin{split}\label{18}
(x+r)(x+r+\lambda )(x+r+2\lambda )+\cdots(x+r+(n-1)\lambda) = \sum_{k=0}^n \genfrac[]{0pt}{1}{n+r}{k+r}_{r,\lambda } x^k.
\end{split}\end{equation}

Note that $\lim_{\lambda \rightarrow 1} \genfrac[]{0pt}{1}{n+r}{k+r}_{r,\lambda } = \genfrac[]{0pt}{1}{n+r}{k+r}_r $, $(n \geq k \geq 0)$.

By Theorem 2.1 and \eqref{18}, we get
\begin{equation}\begin{split}\label{19}
(x-r)_{n,\lambda }=\sum_{k=0}^n S_{1,\lambda }^{(-r)}(n,k)x^k,
\end{split}\end{equation}
and
\begin{equation}\begin{split}\label{20}
(x-r)_{n,\lambda } = \sum_{k=0}^n (-1)^{n-k} \genfrac[]{0pt}{1}{n+r}{k+r}_{r,\lambda }x^k.
\end{split}\end{equation}

From \eqref{19} and \eqref{20}, we can easily derive the following equation \eqref{21}.
\begin{equation}\begin{split}\label{21}
S_{1,\lambda }^{(-r)}(n,k) = (-1)^{n-k} \genfrac[]{0pt}{1}{n+r}{k+r}_{r,\lambda },\,\,(n \geq k \geq 0).
\end{split}\end{equation}

For $n \geq 1$, by Theorem 2.1, we get
\begin{equation}\begin{split}\label{22}
(x+r)_{n+1,\lambda } = \sum_{k=0}^{n+1} S_{1,\lambda }^{(r)}(n+1,k) x^k = \sum_{k=1}^{n+1} S_{1,\lambda }^{(r)}(n+1,k)x^k + (r)_{n+1,\lambda }.
\end{split}\end{equation}

On the other hand, by \eqref{02}, we get
\begin{equation}\begin{split}\label{23}
&(x+r)_{n+1,\lambda }= (x+r)_{n,\lambda }(x+r-n\lambda )\\
&=x \sum_{k=0}^n S_{1,\lambda }^{(r)}(n,k)x^k - (n\lambda -r) \sum_{k=0}^n S_{1,\lambda }^{(r)}(n,k)x^k\\
&=\sum_{k=1}^n S_{1,\lambda }^{(r)}(n,k-1) x^k - \sum_{k=1}^n (n\lambda -r) S_{1,\lambda }^{(r)}(n,k) x^k + (r-n\lambda )(r)_{n,\lambda}+x^{n+1}\\
&= \sum_{k=1}^n \left\{ S_{1,\lambda }^{(r)}(n,k-1) -(n\lambda -r) S_{1,\lambda }^{(r)} (n,k) \right\} x^k + (r)_{n+1,\lambda}+x^{n+1}.
\end{split}\end{equation}

Therefore, by Theorem 2.1 and \eqref{23}, we obtain the following theorem.

\begin{thm}

For $1 \leq k \leq n$, we have
\begin{equation*}\begin{split}
S_{1,\lambda }^{(r)}(n+1,k) = S_{1,\lambda }^{(r)}(n,k-1) -(n\lambda -r) S_{1,\lambda }^{(r)} (n,k).
\end{split}\end{equation*}
\end{thm}

From \eqref{13}, we note that
\begin{equation}\begin{split}\label{24}
\frac{1}{k!} \left( \frac{\log(1+\lambda t)}{\lambda } \right)^k (1+\lambda t)^{\frac{r}{\lambda }} &= \frac{1}{k!} \left( \frac{\log(1+\lambda t)}{\lambda } \right)^k \sum_{l=0}^\infty \frac{r^l}{l!} \left( \frac{\log(1+\lambda t)}{\lambda } \right)^l \\
&= \sum_{l=0}^\infty {k+l \choose l} r^l \frac{1}{(k+l)!} \left( \frac{\log(1+\lambda t)}{\lambda } \right)^{k+l} \\
&= \sum_{l=0}^\infty {k+l \choose l} r^l \sum_{n=k+l}^\infty S_{1,\lambda }(n,k+l) \frac{t^n}{n!}\\
&= \sum_{l=0}^\infty r^l {k+l \choose l}  \sum_{n=l}^\infty S_{1,\lambda }(n+k,k+l) \frac{t^{n+k}}{(n+k)!} \\
&= \sum_{n=0}^\infty \left( \frac{n! t^k}{(n+k)!} \sum_{l=0}^n r^l {k+l \choose l} S_{1,\lambda }(n+k,k+l) \right) \frac{t^n}{n!}.
\end{split}\end{equation}

On the other hand, we have
\begin{equation}\begin{split}\label{25}
\frac{1}{k!} \left( \frac{\log(1+\lambda t)}{\lambda } \right)^k (1+\lambda t)^{\frac{r}{\lambda }} &= \frac{t^k}{k!} \left( \frac{\log(1+\lambda t)}{\lambda t} \right)^k (1+\lambda t)^{\frac{r}{\lambda }} \\
&= \left( \sum_{l=0}^\infty D_l^{(k)} \frac{\lambda^l t^l}{l!} \right) \left( \sum_{m=0}^\infty (r)_{m,\lambda } \frac{t^m}{m!} \right) \frac{t^k}{k!}\\
&=\left( \sum_{n=0}^\infty \sum_{l=0}^n {n \choose l} D_l^{(k)}\lambda^l (r)_{n-l,\lambda } \frac{t^n}{n!} \right) \frac{t^k}{k!}.
\end{split}\end{equation}

Thus, by \eqref{24} and \eqref{25}, we get
\begin{equation}\begin{split}\label{26}
\sum_{l=0}^n r^l \frac{{k+l \choose l}}{{n+k \choose n}} S_{1,\lambda }(n+k,k+l) = \sum_{l=0}^n {n \choose l}D_l^{(k)}\lambda^l(r)_{n-l,\lambda }.
\end{split}\end{equation}

Therefore, by \eqref{26}, we obtain the following theorem.

\begin{thm}
For $n \geq 0$, we have
\begin{equation*}\begin{split}
\sum_{l=0}^n {n \choose l}D_l^{(k)} \lambda^l(r)_{n-l,\lambda } = \sum_{l=0}^n  \frac{{k+l \choose l}}{{n+k \choose n}}r^l S_{1,\lambda }(n+k,k+l).
\end{split}\end{equation*}
\end{thm}

Now, we observe that
\begin{equation}\begin{split}\label{27}
\frac{1}{k!} \left( \frac{\log(1+\lambda t)}{\lambda } \right)^k (1+\lambda t)^{\frac{r}{\lambda }} &= \left(\sum_{l=0}^\infty (r)_{l,\lambda}\frac{t^l}{l!}\right)\frac{1}{k!} \left( \frac{\log(1+\lambda t)}{\lambda } \right)^k  \\
&= \sum_{n=k}^\infty \left( \sum_{m=k}^n {n \choose m} S_{1,\lambda }(m,k) (r)_{n-m,\lambda } \right) \frac{t^n}{n!}.
\end{split}\end{equation}

Therefore, by \eqref{13} and \eqref{27}, we obtain the following theorem.

\begin{thm}
For $n,k \geq 0$, with $n \geq k$, we have
\begin{equation*}\begin{split}
S_{1,\lambda }^{(r)}(n,k) = \sum_{m=k}^n {n \choose m}  (r)_{n-m,\lambda }  S_{1,\lambda }(m,k).
\end{split}\end{equation*}
\end{thm}

From \eqref{13}, we note that
\begin{equation}\begin{split}\label{28}
&\frac{1}{m!} \left( \frac{\log(1+\lambda t)}{\lambda } \right)^m
\frac{1}{k!} \left( \frac{\log(1+\lambda t)}{\lambda } \right)^k (1+\lambda t)^{\frac{r}{\lambda }} \\
&= \frac{(m+k)!}{m!k!} \frac{1}{(m+k)!} \left( \frac{\log(1+\lambda t)}{\lambda } \right)^{m+k} (1+\lambda t)^{\frac{r}{\lambda }} \\
&={m+k \choose m} \sum_{n=m+k}^\infty S_{1,\lambda }^{(r)}(n,m+k) \frac{t^n}{n!}.
\end{split}\end{equation}

On the other hand,
\begin{equation}\begin{split}\label{29}
&\frac{1}{m!} \left( \frac{\log(1+\lambda t)}{\lambda } \right)^m
\frac{1}{k!} \left( \frac{\log(1+\lambda t)}{\lambda } \right)^k (1+\lambda t)^{\frac{r}{\lambda }} \\
&= \left( \sum_{l=m}^\infty S_{1,\lambda }(l,m) \frac{t^l}{l!} \right) \left( \sum_{j=k}^\infty S_{1,\lambda }^{(r)}(j,k) \frac{t^j}{j!} \right)\\
&= \sum_{n=m+k}^\infty \left( \sum_{l=k}^{n-m}{n \choose l} S_{1,\lambda }^{(r)}(l,k) S_{1,\lambda }(n-l,m) \right) \frac{t^n}{n!}.
\end{split}\end{equation}

Therefore, by \eqref{28} and \eqref{29}, we obtain the following theorem.

\begin{thm}
For $m,n,k \geq 0$ with $n \geq m+k$, we have
\begin{equation*}\begin{split}
{m+k \choose m} S_{1,\lambda }^{(r)}(n,m+k) = \sum_{l=k}^{n-m}{n \choose l} S_{1,\lambda }(l,k) S_{1,\lambda }(n-l,m).
\end{split}\end{equation*}
\end{thm}

By \eqref{12}, we get
\begin{equation}\begin{split}\label{30}
\sum_{n=k}^\infty S_{1,\lambda }(n,k) \frac{t^n}{n!} &=
\frac{1}{k!}\left( \frac{\log(1+\lambda t)}{\lambda } \right)^k (1+\lambda t)^{\frac{r}{\lambda }} (1+\lambda t)^{-\frac{r}{\lambda }}\\
&=\left( \sum_{l=k}^\infty S_{1,\lambda }^{(r)}(l,k) \frac{t^l}{l!} \right) \left( \sum_{m=0}^\infty { -\frac{r}{\lambda } \choose m} \lambda ^m t^m \right) \\
&=\left( \sum_{l=k}^\infty S_{1,\lambda }^{(r)}(l,k) \frac{t^l}{l!} \right) \left( \sum_{m=0}^\infty (-1)^m (r+(m-1)\lambda )_{m,\lambda} \frac{t^m}{m!} \right) \\
&=\sum_{n=k}^\infty \left( \sum_{l=k}^n {n \choose l} S_{1,\lambda }^{(r)}(l,k) (-1)^{n-l} (r+(n-l-1)\lambda)_{n-l,\lambda } \right) \frac{t^n}{n!}.
\end{split}\end{equation}

Comparing the coefficients on both sides of \eqref{30}, we have the following theorem.

\begin{thm}
For $n,k \geq 0$, with $n \geq k$, we have
\begin{equation*}\begin{split}
S_{1,\lambda }(n,k) = \sum_{l=k}^n {n \choose l} S_{1,\lambda }^{(r)}(l,k) (-1)^{n-l} (r+\lambda (n-l-1))_{n-l,\lambda }.
\end{split}\end{equation*}
\end{thm}

From \eqref{09}, we have
\begin{equation}\begin{split}\label{32}
&\frac{1}{k!}\left( \frac{\log(1+\lambda t)}{\lambda } \right)^k (1+\lambda t)^{\frac{r}{\lambda }} = \frac{t^k}{k!} \left( \frac{\log(1+\lambda t)}{\lambda t} \right)^k (1+\lambda t)^{\frac{r}{\lambda }}\\
&= \frac{t^k}{k!} \left( \sum_{m=0}^\infty D_m^{(k)}\lambda ^m \frac{t^m}{m!} \right) \left( \sum_{l=0}^\infty (r)_{l,\lambda }\frac{t^l}{l!} \right)\\
&= \frac{t^k}{k!} \sum_{n=0}^\infty  \left( \sum_{m=0}^n {n \choose m} D_m^{(k)} \lambda ^m (r)_{n-m,\lambda } \right)
 \frac{t^n}{n!}.
\end{split}\end{equation}

On the other hand, by \eqref{13}, we get
\begin{equation}\begin{split}\label{33}
\frac{1}{k!}\left( \frac{\log(1+\lambda t)}{\lambda } \right)^k (1+\lambda t)^{\frac{r}{\lambda }} &= \sum_{n=k}^\infty S_{1,\lambda }^{(r)}(n,k) \frac{t^n}{n!}\\
& = \frac{t^k}{k!} \sum_{n=0}^\infty S_{1,\lambda }^{(r)}(n+k,k) \frac{n!k!}{(n+k)!} \frac{t^n}{n!}.
\end{split}\end{equation}

Thus, by comparing the coefficients on both sides of \eqref{32} and \eqref{33}, we get
\begin{equation}\begin{split}\label{34}
\sum_{m=0}^n {n \choose m}D_m^{(k)}\lambda ^m (r)_{n-m,\lambda } = \frac{1}{{n+k \choose n} }S_{1,\lambda }^{(r)}(n+k,k).
\end{split}\end{equation}

Therefore, by \eqref{34}, we obtain the following theorem.
\begin{thm}
For $n,k \geq 0$, we have
\begin{equation*}\begin{split}
S_{1,\lambda }^{(r)}(n+k,k) = {n+k \choose n} \sum_{m=0}^n {n \choose m}D_m^{(k)}\lambda^m (r)_{n-m,\lambda }.
\end{split}\end{equation*}
\end{thm}

From \eqref{09}, we note that
\begin{equation}\begin{split}\label{35}
\frac{1}{k!}\left( \frac{\log(1+\lambda t)}{\lambda } \right)^k (1+\lambda t)^{\frac{r}{\lambda }}&= \frac{t^k}{k!} \left( \frac{\log(1+\lambda t)}{\lambda t} \right)^k (1+\lambda t)^{\frac{r}{\lambda }}\\
&= \frac{t^k}{k!} \sum_{n=0}^\infty \lambda ^n D_n^{(k)} (\tfrac{r}{\lambda }) \frac{t^n}{n!}.
\end{split}\end{equation}

By \eqref{33} and \eqref{35}, we get
\begin{equation}\begin{split}\label{36}
S_{1,\lambda }^{(r)}(n+k,k) = \lambda ^n \frac{(n+k)!}{n!k!} D_n^{(k)} (\tfrac{r}{\lambda }) = \lambda ^n {n+k \choose n} D_n^{(k)} (\tfrac{r}{\lambda }),\,\,(n \geq 0).
\end{split}\end{equation}

In particular, for $r=0$, from \eqref{30} and \eqref{36} we have

\begin{equation}\begin{split}\label{37}
&\lambda ^n {n+k \choose k} D_n^{(k)} = S_{1,\lambda }(n+k,k)\\
&= \sum_{l=k}^{n+k} {n+k \choose l} S_{1,\lambda }^{(r)}(l,k) (-1)^{n+k-l} (r+(n+k-l-1)\lambda)_{n+k-l,\lambda },
\end{split}\end{equation}
where $n,k \ge 0$.

Therefore, by \eqref{37}, we obtain the following theorem.
\begin{thm}
For $n,k \geq 0$, we have
\begin{equation*}\begin{split}
\lambda ^n {n+k \choose k} D_n^{(k)}=\sum_{l=k}^{n+k} {n+k \choose l} S_{1,\lambda }^{(r)}(l,k) (-1)^{n+k-l} (r+ (n+k-l-1)\lambda)_{n+k-l,\lambda }.
\end{split}\end{equation*}
In addition,
\begin{equation*}\begin{split}
& D_n^{(k)}=\frac{1}{{n+k \choose k}} \sum_{l=k}^{n+k} {n+k \choose l} {l \choose k} \left( \frac{1}{\lambda } \right)^{n+k-l}\\
&\times (r+ (n+k-l-1)\lambda )_{n+k-l,\lambda} (-1)^{n+k-l} D_{l-k}^{(k)}(\tfrac{r}{\lambda }).
\end{split}\end{equation*}
\end{thm}

Now, we observe that
\begin{equation}\begin{split}\label{38}
\sum_{n=k}^\infty S_{1,\lambda }(n,k) \frac{t^n}{n!} &=
\frac{1}{k!} \left( \frac{\log(1+\lambda t)}{\lambda } \right)^k (1+\lambda t)^{\frac{r}{\lambda }} e^{-\frac{r}{\lambda }\log(1+\lambda t)}\\
&=\left(\sum_{l=k}^\infty S_{1,\lambda }^{(r)}(l,k) \frac{t^l}{l!}\right) \sum_{m=0}^\infty (-1)^m r^m \frac{1}{m!} \left( \frac{\log(1+\lambda t)}{\lambda } \right)^m \\
&=\left( \sum_{l=k}^\infty S_{1,\lambda }^{(r)}(l,k) \frac{t^l}{l!} \right) \left( \sum_{m=0}^\infty (-1)^m r^m \sum_{j=m}^\infty S_{1,\lambda }(j,m) \frac{t^j}{j!} \right) \\
&=\sum_{n=k}^\infty \left(
\sum_{j=0}^{n-k} \sum_{m=0}^j {n \choose j}(-1)^m r^m S_{1,\lambda }(j,m) S_{1,\lambda }(n-j,k)  \right) \frac{t^n}{n!}.
\end{split}\end{equation}

Therefore, by comparing the coefficients on both sides of \eqref{38}, we obtain the following theorem
\begin{thm}
For $n,k \geq 0$, with $n \geq k$, we have
\begin{equation*}\begin{split}
S_{1,\lambda }(n,k) = \sum_{j=0}^{n-k} \sum_{m=0}^j {n \choose j}(-1)^m r^m S_{1,\lambda }(j,m) S_{1,\lambda }(n-j,k)  .
\end{split}\end{equation*}
\end{thm}

For $m,n \ge 0$, we define $\lambda $-analogues of the Whitney's type $r$-Stirling numbers of the first kind as
\begin{equation}\begin{split}\label{39}
(mx+r)_{n,\lambda } &= (mx+r)(mx+r-\lambda )(mx+r-2\lambda )\cdots(mx+r-(n-1)\lambda )\\
&= \sum_{k=0}^n T_{1,\lambda }^{(r)}(n,k|m) x^k.
\end{split}\end{equation}

By \eqref{39}, we get
\begin{equation}\begin{split}\label{40}
\sum_{n=0}^\infty (mx+r)_{n,\lambda } \frac{t^n}{n!} &= \sum_{n=0}^\infty \left( \sum_{k=0}^n T_{1,\lambda }^{(r)}(n,k|m) x^k \right) \frac{t^n}{n!}\\
&= \sum_{k=0}^\infty \left(\sum_{n=k}^\infty T_{1,\lambda }^{(r)}(n,k|m) \frac{t^n}{n!} \right) x^k.
\end{split}\end{equation}

On the other hand, by binomial expansion, we get
\begin{equation}\begin{split}\label{41}
\sum_{n=0}^\infty (mx+r)_{n,\lambda }\frac{t^n}{n!}&= \sum_{n=0}^\infty {mx+r \choose n}_\lambda t^n \\
&= (1+\lambda t)^{\frac{mx+r}{\lambda }} = (1+\lambda t)^{\frac{r}{\lambda }} e^{mx(\frac{\log(1+\lambda t)}{\lambda })}\\
& = \sum_{k=0}^\infty \frac{m^k}{k!} \left( \frac{\log(1+\lambda t)}{\lambda } \right)^k (1+\lambda t)^{\frac{r}{\lambda }} x^k.
\end{split}\end{equation}

Comparing the coefficients on both sides of \eqref{40} and \eqref{41}, the generating function for $T_{1,\lambda }^{(r)}(n,k|m)$, $(n,k \geq 0)$, is given by
\begin{equation}\begin{split}\label{42}
\frac{m^k}{k!} \left( \frac{\log(1+\lambda t)}{\lambda } \right)^k (1+\lambda t)^{\frac{r}{\lambda }} = \sum_{n=k}^\infty T_{1,\lambda }^{(r)}(n,k|m) \frac{t^n}{n!}.
\end{split}\end{equation}

From \eqref{13} and \eqref{42}, we note that
\begin{equation}\begin{split}\label{43}
S_{1,\lambda }^{(r)}(n,k) = \frac{1}{m^k} T_{1,\lambda }^{(r)}(n,k|m),\,\,(n \geq k \geq 0).
\end{split}\end{equation}

It is known that the $r$-Whitney numbers are defined as
\begin{equation}\begin{split}\label{44}
(mx+r)^n = \sum_{k=0}^n m^k W_{m,r}(n,k) (x)_k, \quad (\text{see}\,\, [3]).
\end{split}\end{equation}

By \eqref{03}, we get

\begin{equation}\begin{split}\label{45}
(mx+r)_{n,\lambda }&= \sum_{l=0}^n S_{1,\lambda }(n,l) (mx+r)^l \\
&= \sum_{l=0}^n S_{1,\lambda }(n,l) \sum_{j=0}^l m^j W_{m,r}(l,j)(x)_j\\
&=\sum_{j=0}^n \sum_{l=j}^n S_{1,\lambda }(n,l) m^j W_{m,.r} (l,j) (x)_j\\
&= \sum_{j=0}^n \sum_{l=j}^n S_{1,\lambda }(n,l) m^j W_{m,r}(l,j) \sum_{k=0}^j S_1 (j,k) x^k\\
&=\sum_{ k=0}^n \left( \sum_{ j=k}^n \sum_{l=j}^n S_{1,\lambda }(n,l) S_1(j,k) m^j W_{m,r}(l,j) \right) x^k.
\end{split}\end{equation}

Therefore, by \eqref{39} and \eqref{45}, we obtain the following theorem.

\begin{thm}
For $n,k \geq 0$, with $n \geq k$, we have
\begin{equation*}\begin{split}
T_{1,\lambda }^{(r)}(n,k|m) = \sum_{ j=k}^n \sum_{l=j}^n S_{1,\lambda }(n,l) S_1(j,k) m^j W_{m,r}(l,j) .
\end{split}\end{equation*}
\end{thm}

\end{document}